\documentclass[manualthm,nocrop,noinfo,final]{pem_v1}

\journal{PEM}
\firstpage{1}
\art-id{PEMS 220731-Burkhart}
\slugauthor{xxx}
\slugrefno{XXX}

\theoremstyle{headfive}

\newtheorem{theorem}{Theorem}

\newtheorem{corollary}{Corollary}
\newtheorem*{schur-zassenhaus}{Theorem}
\newtheorem*{higman}{Corollary}
\newtheorem*{glauberman}{Lemma}

\theoremstyle{headsix}

\newcommand{\semi}{\rtimes}
\newcommand{\iso}{\cong}
\newcommand{\normal}{\vartriangleleft}
\renewcommand{\phi}{\varphi}
\providecommand{\abs}[1]{\ensuremath{\left\lvert#1\right\rvert}}

\DeclareMathOperator{\Fix}{Fix}

\begin{document}

\title[Conjugacy conditions for supersoluble complements]
{Conjugacy conditions for supersoluble complements of an abelian base and a
	fixed point result for non-coprime actions}

\author[M.\,C. Burkhart]{
\FMAUGRP{
\FMAUT{
\FMFNM{Michael C.}\FMSEP{ }
\FMSNM{Burkhart}
}}}

\affiliation{
\FMAFFGRP{
\FMAFF{
\FMINST{University of Cambridge},
\FMCYNAME{Cambridge},
\FMCOUNTRY{United Kingdom}
\email{mcb93@cam.ac.uk}
}}}

\date{Received 31 July 2022}

\begin{abstract}
	We demonstrate that two supersoluble complements of an abelian base in a finite
	split extension are conjugate if and only if, for each prime $p$, a Sylow
	$p$-subgroup of one complement is conjugate to a Sylow $p$-subgroup of the
	other. As a corollary, we find that any two supersoluble complements of an
	abelian subgroup $N$ in a finite split extension $G$ are conjugate if and only
	if, for each prime $p$, there exists a Sylow $p$-subgroup $S$ of $G$ such that
	any two complements of $S\cap N$ in $S$ are conjugate in $G$. In particular,
	restricting to supersoluble groups allows us to ease D.\,G. Higman's stipulation that
	the complements of $S\cap N$ in $S$ be conjugate within $S$. We then consider
	group actions and obtain a fixed point result for non-coprime actions analogous
	to Glauberman's lemma.
\end{abstract}


\classification[{\rm 2020} Mathematics subject classification]{Primary 20E22; 20F16; 20E45 \break Secondary 05E18}

\maketitle

\section{Introduction}
\label{s:background}

A finite group $G$ splits over a normal subgroup $N$ if there exists a
complement $H$ such that $G\iso N\semi H$ forms a semidirect product. When
$\abs{N}$ and $[G:N]$ are coprime, Schur proved a complement must exist and
Zassenhaus demonstrated, under the additional assumption\footnote{rendered
	unnecessary by Feit and Thompson's theorem that every finite group of odd order
	is soluble} that either $N$ or $G/N$ is soluble, that such a complement is
unique up to conjugacy. When $N$ is abelian, Gasch\"utz showed a complement
exists if and only if, for each prime $p$, there exists a Sylow $p$-subgroup
$S$ of $G$ that splits over $S\cap N$. Higman considered the question of
conjugacy for such complements and found that~\cite[cor. 2]{Hig54}
\begin{higman}[D.\,G. Higman]
	Let $G$ be a split extension of an abelian subgroup $N$. If for each prime $p$
	there is a Sylow $p$-subgroup $S$ of $G$ such that any two complements of
	$S\cap N$ in $S$ are conjugate in $S$, then any two complements of $N$ in $G$
	are conjugate.
\end{higman}

In this note, we establish similar criteria for complements of an abelian base
to be conjugate, loosening the restrictiveness of the conjugacy condition for
the Sylow subgroups. To this end, we consider supersoluble groups, namely those
groups $G$ that possess a series $1=G_0 \normal G_1 \normal \dotsb \normal G_{s-1}
	\normal G_s = G$ with each factor $G_i$ normal in $G$ and each quotient
$G_{i+1}/G_i$ cyclic, for $0\leq i < s$. All finite nilpotent groups are
supersoluble and in turn all supersoluble groups are soluble. We find that%
\begin{theorem}
	\label{thm}
	Suppose the supersoluble subgroups $H$ and $H'$ each complement a normal
	abelian subgroup $N$ in a finite group $G$. If for each prime $p$, $H'$
	contains a conjugate of some Sylow $p$-subgroup of $H$, then $H$ and $H'$ are
	conjugate.
\end{theorem}
We then offer the following corollary, more in the spirit of Higman's work,
\begin{corollary}
	\label{cor_higman}
	Let $G$ be a split extension of an abelian subgroup $N$ such that $G/N$ is supersoluble. If for
	each prime $p$ there is a Sylow $p$-subgroup $S$ of $G$ such that any two
	complements of $S\cap N$ in $S$ are conjugate in $G$, then any two complements
	of $N$ in $G$ are conjugate.
\end{corollary}
If $G$ is supersoluble, then $G/N$ is as well. Thus, for supersoluble groups, complements of $S\cap N$ in $S$ need not be
conjugate strictly within $S$ but merely in the full group $G$.

We then consider the fixed points of groups with operator groups, where
Glauberman's result on the existence of fixed points for coprime
actions~\cite[thm. 4]{Gla64} states
\begin{glauberman}[Glauberman]
	Given a finite group $H$ acting via automorphisms on a finite group $N$,
	suppose the induced semidirect product $N\semi H$ acts on some non-empty set
	$\Omega$ where the action of $N$ is transitive. If the orders of $H$ and $N$
	are coprime, then there exists an $H$-invariant element $\omega\in\Omega$.
\end{glauberman}
We offer an analogous fixed point result for non-coprime actions. In
particular, theorem~\ref{thm} and Gasch\"utz's existence result allow us to
demonstrate
\begin{corollary}
	\label{cor_glauberman}
	Given a finite supersoluble group $H$ acting via automorphisms on a finite
	abelian group $N$, suppose the induced semidirect product $N\semi H$ acts on
	some non-empty set $\Omega$ where the action of $N$ is transitive. If for each
	prime $p$, a Sylow $p$-subgroup of $H$ fixes some element of $\Omega$, then
	there exists an $H$-invariant element $\omega\in\Omega$.
\end{corollary}
While Isaacs previously considered non-coprime actions and outlined some
conditions for the existence of a fixed point~\cite{Isa68}, results of this
nature appear to be rare in the literature.

\subsection{Outline}
We proceed as follows. In the remainder of this section, we introduce some
notation before proving the main theorem in section~\ref{s:pf_of_thm}. We then
prove the corollaries in section~\ref{s:pf_of_cors} and make concluding remarks
in section~\ref{s:conclusion}.

\subsection{Notation}
All groups in this note are assumed to be finite. We use standard notation
readily found in any introductory textbook on group theory, such as
Isaacs'~\cite{Isa08}. We let groups act from the left and denote conjugation by
$g^\gamma = \gamma^{-1}g\gamma$ for $g,\gamma\in G$. For subgroups $J$ and $K$
of $G$, we let $N_K(J)$ and $C_K(J)$ denote the elements of $K$ that normalise
and centralise $J$, respectively. In the semidirect product $G \iso N \semi H$,
we note that any two complements of $N$ are necessarily isomorphic (to the
quotient group $G/N$) and that $N_K(J) = C_K(J)$ whenever both $K\leq N$ and $J \leq
	H$, as then $[J, N_K(J)] \leq J \cap N = 1$.

\section{Proof of main theorem}
\label{s:pf_of_thm}

We first restrict $N$ to be an abelian $p$-subgroup and then consider the
general case that $N$ is abelian. In both cases that follow, we assume we are
given supersoluble complements $H$ and $H'$ of a nontrivial subgroup $N\normal
	G$ satisfying the hypotheses of theorem~\ref{thm}. The theorem would hold
trivially if $H$ were a $q$-subgroup for some prime $q$, so without loss, at
least two distinct primes divide $\abs{H}$. Zassenhaus's result further implies
$p$ must be one of these primes. As any two Sylow subgroups of $H$ are
conjugate within $H$ for a given prime, our assumptions imply that $H'$
contains some conjugate of every Sylow subgroup of $H$ where each such
conjugate is itself a Sylow subgroup of $H'$. Furthermore, as any element $g\in
	G$ can be uniquely written $g=hn$ for some $n\in N$ and $h\in H$, we may assume
that for each prime $p$, there exists a Sylow $p$-subgroup $P$ of $H$ such that
$P^n\leq H'$ for some $n\in N$.

\subsection{Case that $N$ is a $p$-subgroup}
	We first suppose $N$ is a $p$-subgroup for some prime $p$ and induct on the
	order of $G$. As $H$ is supersoluble, we may find a normal subgroup $H_0
		\normal H$ of index $q$~\cite[exer. 3B.9]{Isa08}, where $q$ is the
	smallest prime divisor of $\abs{H}$. Under the inductive hypothesis, $H_0$ is
	conjugate to some subgroup of $H'$ and so the left action of $H_0$ on the
	cosets $\Omega=G/H'$ has a fixed point, say $n'H'$. If $H$ itself does not fix any
	point of $\Omega$, then the nonempty set $\Fix_{\Omega}(H_0)$ of points in
	$\Omega$ fixed by $H_0$ may be partitioned into orbits of $H$, each having
	cardinality $[H:H_0] = q$. On the other hand, the map from $N_N(H_0)$ to
	$\Fix_{\Omega}(H_0)$ given by $n\mapsto n\cdot n'H'$ provides a bijective
	correspondence between these two sets, and as $N_N(H_0) \leq N$ is a $p$-group,
	we may assume without loss that $q=p$.  In particular, $H_0$ contains a
	$p'$-Hall subgroup of $H$, say $M$, such that $C_N(M) \supseteq C_N(H_0) = N_N(H_0)$ is
	nontrivial. Furthermore, we may assume without loss that $M\normal H$~\cite[exer. 3B.8]{Isa08} so 
	$N_0 := [N, M]  \normal G$ where $N=N_0 \times C_N(M)$ by a theorem of
	Fitting~\cite[thm. 4.34]{Isa08}.  Consequently, $N_0$ is a strict subgroup of $N$.
	
	Let $P$ be a Sylow $p$-subgroup of $H$. By hypothesis, there exists some 
	$n \in N$ such that $P^n \leq H'$. Replacing $H$ by $H^n$ in the statement 
	of the theorem if necessary, we may assume without loss that $n$ is trivial,
	 i.e. that $P \leq H'$.

	Glauberman's lemma implies the left action of $NM/N_0$ on the
	cosets $G/N_0H'$ has an $M$-invariant element. As $N=N_0 \times C_N(M)$, 
	it follows that $\overline{N} = \overline{C_N(M)}$ in $\overline
		G = G/N_0$. Furthermore, as fixed points come from fixed points in coprime
	actions~\cite[cor. 3.28]{Isa08}, $\overline{N} = C_{\overline
				N}(M)$ so that $M$ fixes every element of $G/N_0 H'$. In particular, $M\leq N_0H'$
	so that $N_0M$ acts on the cosets $N_0H'/H'$. Glauberman's lemma implies this
	second action also has an $M$-invariant element so that $M^{n_0} \leq H'$ for
	some $n_0 \in N_0$. Thus, as $H\leq N_0H'$, we may apply the inductive hypothesis
	in $N_0H'$ to conclude that $H$ and $H'$ are conjugate.

\subsection{General case that $N$ is abelian}
	We again induct on the order of $G$. Building upon the results of our first
	case, we now assume $\abs{N}$ has multiple prime divisors. In particular, we
	may write $N=N_1\times N_2$ where $N_1$ and $N_2$ are characteristic in $N$ and
	thus normal in $G$. (We may, for example, let $N_1$ be the first primary
	component of $N$ and $N_2$ be the product of the other primary components.) In
	$G/N_1$, the inductive hypothesis implies the action of $H$ on the left
	cosets $G/N_1H'$ has a fixed point, say $n_2N_1H'$ for some $n_2\in N_2$.
	Analogously, in $G/N_2$, we find that the action of $H$ on $G/N_2H'$ must also have a fixed point, say
	$n_1N_2H'$ for some $n_1\in N_1$. Consequently, we conclude that the action of
	$H$ on $G/H'$ fixes $n_1n_2H'$ so that $H$ and $H'$ are conjugate.

\section{Proofs of corollaries}
\label{s:pf_of_cors}

In this section, we provide proofs for the corollaries stated in
section~\ref{s:background}.

\subsection{Proof of Corollary~\ref{cor_higman}}
	Given a group $G$ satisfying the hypotheses of corollary~\ref{cor_higman}, we
	suppose by way of contradiction that $H$ and $H'$ complement $N$ in $G$ but
	fail to be conjugate. By theorem~\ref{thm}, there exists a prime $p$ and Sylow
	$p$-subgroups $P$ and $P'$ of $H$ and $H'$, respectively, such that $P$ and
	$P'$ fail to be conjugate in $G$. Let $L$ be the unique Sylow $p$-subgroup of
	$N$. Then $LP$ and $LP'$ are Sylow $p$-subgroups of $G$ and thus are conjugate
	to the Sylow $p$-subgroup $S$ of $G$ specified in the statement of the
	corollary. In particular, $LP^{g_1}=L(P')^{g_2}=S$ for some $g_1,g_2\in G$
	where $S\cap N = L$, so that by assumption, $P^{g_1}$ and $({P'})^{g_2}$ are
	conjugate in $G$, a contradiction.

\subsection{Proof of Corollary~\ref{cor_glauberman}}
	Suppose $G=N \semi H$ acts on a set $\Omega$ according to the hypotheses of
	corollary~\ref{cor_glauberman}. Fix some $\alpha\in\Omega$ and consider the
	point stabilizer $G_\alpha$ of $\alpha$ in $G$. As $N$ is transitive on
	$\Omega$, for any $g\in G$, we have $g\cdot\alpha = n\cdot \alpha$ for some
	$n\in N$, so that $n^{-1}g \in G_\alpha$. Thus, $G=NG_\alpha$. As $N\normal G$,
	it follows that $G_\alpha\cap N \normal G_\alpha$.

	We claim $G_\alpha$ splits over $G_\alpha\cap N$. By Gasch\"utz's theorem, it
	suffices to show that for each prime $p$, there exists a Sylow p-subgroup S of
	$G_\alpha$ that splits over $S\cap N$. Fix $p$, an arbitrary prime. By
	hypothesis, there exists some $n\in N$ and a Sylow $p$-subgroup $P$ of $H$ such
	that $P^n \leq G_\alpha$. Let $L$ be the Sylow $p$-subgroup of $G_\alpha\cap N$.
	As $\abs{G_\alpha} = \abs{G_\alpha\cap N} [G:N]$, it follows that $S=LP^n$ is a
	Sylow subgroup of $G_\alpha$ and $P^n$ is a complement of $S \cap N = L$ in
	that subgroup.

	Thus $G_\alpha$ splits and we may let $H'$ complement $G_\alpha\cap N$ in
	$G_\alpha$. It follows that $\abs{H'} = [G_\alpha: G_\alpha\cap N] = [G:N]$ so
	$H'$ complements $N$ in $G$. Theorem~\ref{thm} then implies that $H'=H^n$ for
	some $n\in N$ so that $H$ fixes $n \cdot \alpha$.

\section{Concluding remarks}
\label{s:conclusion}

We conclude with some historical context.

\subsection{Conjugacy of complements}
D.\,G. Higman's work included another result on the conjugacy of complements
that we mention here for completeness. In a split extension $G\iso N\semi H$ of
an abelian subgroup $N$, Higman considered an intermediate subgroup $N\leq B
	\leq G$ of index $[G:B]=b$ such that the map $n\mapsto n^b$ for $n\in N$ has an
inverse. (Such an inverse always exists if $b$ and $\abs{N}$ are coprime, for
example.) He showed that if any two complements of $N$ in $B$ are conjugate in
$B$, then any two complements of $N$ in $G$ are conjugate in $G$. Furthermore,
given two specific complements $H$ and $H'$ of $N$ in $G$, he showed that $H$
and $H'$ are conjugate in $G$ if and only if $H\cap B$ and $H'\cap B$ are
conjugate in $B$~\cite[cor. 1]{Hig54}.

\subsection{Fixed points}
Glauberman's original paper~\cite{Gla64} provided two conditions for
$N$-transitive actions of split extensions $G\iso N\semi H$ that jointly imply
the existence of an $H$-invariant element, namely (S) if $W\leq G$ satisfies
$WN=G$ then $W$ splits over $W\cap N$, and (Z) any two complements of $N$ in
$G$ are conjugate.\footnote{conditions (S) and (Z) correspond to Schur's
	existence and Zassenhaus's conjugacy results for the complements of normal Hall
	subgroups}~ The hypotheses of our corollary~\ref{cor_glauberman} do not
necessarily imply (Z) but rather only that two specific complements of $N$ are
conjugate, namely $H$ and a complement of $G_\alpha$ in $G_\alpha\cap N$, where
$\alpha\in\Omega$ is a member of the set on which $G$ acts. For this reason, we
did not directly apply Glauberman's result but instead mirrored his proof.

To further illustrate the connection between the conjugacy of complements and
the fixed points of certain group actions, we note that
corollary~\ref{cor_glauberman} implies theorem~\ref{thm}. To see this
equivalence, we may apply corollary~\ref{cor_glauberman} to the action
of $N\semi H$ on the set $\Omega=G/H'$ of left cosets for any second complement
$H'$ satisfying the hypotheses of theorem~\ref{thm}.

Finally, we note that in the language of group cohomology,
corollary~\ref{cor_higman} provides a vanishing condition for the first
cohomology group $H^1(H,N)$ of the split extension $G\iso N\semi H$ when $N$ is
abelian and $H$ is supersoluble.

\ack{	We thank Elizabeth Crites, Michael Bate, Charles Eaton, Lenny Taelman, Martin
	Liebeck, Laura Ciobanu, and anonymous reviewers for their thoughtful and insightful feedback.}
	
\compete{The author declares none.}

\end{document}